\documentclass{article}

\usepackage{amssymb, amsmath, latexsym, mathrsfs, graphicx, stmaryrd, amsthm}
\usepackage{url}

\newtheorem{thm}{Theorem}[section]

\newtheorem{exm}[thm]{Example}

\theoremstyle{definition}

\newcommand{\N}{\mathbb{N}}
\newcommand{\Z}{\mathbb{Z}}

\newcommand{\R}{\mathbb{R}}

\newcommand{\LL}{\mathcal D}

\makeatletter
\def\revddots{\mathinner{\mkern1mu\raise\p@
\vbox{\kern7\p@\hbox{.}}\mkern2mu
\raise4\p@\hbox{.}\mkern2mu\raise7\p@\hbox{.}\mkern1mu}}
\makeatother

\renewcommand{\mod}[1]{\,(\text{mod }#1)}

\begin{document}

\title{The Minimum Period of the Ehrhart Quasi-polynomial of a Rational Polytope}

\author{Tyrrell B. McAllister\thanks{Partially supported by NSF VIGRE Grant No. DMS-0135345.} \and Kevin M. Woods\thanks{Partially supported by an NSF Graduate Research Fellowship.}}

\maketitle
\begin{abstract}
If $P\subset \R^d$ is a rational polytope, then $i_P(n):=\#(nP\cap \Z^d)$ is a quasi-polynomial in $n$, called the Ehrhart quasi-polynomial of $P$.  The period of $i_P(n)$ must divide $\LL(P)= \min \{n \in \Z_{> 0} \colon  nP \text{ is an integral polytope}\}$.  Few examples are known where the period is not exactly $\LL(P)$.  We show that for any $\LL$, there is a 2-dimensional triangle $P$ such that $\LL(P)=\LL$ but such that the period of $i_P(n)$ is 1, that is, $i_P(n)$ is a polynomial in $n$.  We also characterize all polygons $P$ such that $i_P(n)$ is a polynomial.  In addition, we provide a counterexample to a conjecture by T. Zaslavsky about the periods of the coefficients of the Ehrhart quasi-polynomial.
\end{abstract}  

\section{Introduction}

An \emph{integral} (respectively, \emph{rational}) \emph{polytope} is a polytope whose vertices have integral (respectively, rational) coordinates.  Given a rational polytope $P \subset \R^d$, the \emph{denominator} of $P$ is
\[
\LL(P) = \min \{n \in \Z_{> 0} \colon  nP \text{ is an integral polytope}\}.
\]

Ehrhart proved (\cite{ehrhart}) that if $P \subset \R^d$ is a rational polytope, then there is a quasi-polynomial function $i_P \colon \Z \mapsto \Z$ with period $\LL(P)$ such that, for $n \geq 0$,
\[
i_P(n) = \#\left( n P \cap \Z^d \right).
\]
In other words, there exist polynomial functions $f_1,\dotsc,f_{\LL(P)}$ such that $i_P(n) = f_j(n)$ for $n \equiv j \mod{\LL(P)}$.  In particular, if $P$ is integral, then $\LL(P) = 1$, so $i_P$ is a polynomial function.

We call $i_P$ the \emph{Ehrhart quasi-polynomial of $P$}.  This counting function satisfies several important properties:
\begin{enumerate}
\item The degree of each $f_j$ is the dimension of $P$.
\item The coefficient of the leading term of each $f_j$ is the volume of $P$, normalized with respect to the sublattice of $\Z^d$ which is the intersection of $\Z^d$ with the affine hull of $P$ (in particular, if $P$ is full dimensional, the coefficient is simply the Euclidean volume of $P$).
\item (Law of Reciprocity) For $n \geq 1$, let
\[
i^\circ_P(n)=\#\left(interior(nP)\cap \Z^d\right).
\]
Then $i^\circ_P(n)=(-1)^d i_P(-n)$.
\end{enumerate}
Properties (1) and (2) were proved by Ehrhart in \cite{ehrhart}.  Property (3) was conjectured by Ehrhart and proved in full generality by I.G. MacDonald in \cite{macdonald}.  For an excellent introduction to Ehrhart quasi-polynomials that includes proofs of all these properties, see \cite{stanley}.

We know that $\LL(P)$ is \emph{a} period of the Ehrhart quasi-polynomial of $P$, but what is the \emph{minimum} period?  Of course, it must divide $\LL(P)$, and it very often equals $\LL(P)$.  Though this is not always the case, very few counterexamples were previously known.   R.P. Stanley (\cite{stanley}, Example 4.6.27) provided an example of a polytope $P$ with denominator $\LL(P) = 2$ where the minimum period is 1, that is, where the Ehrhart quasi-polynomial is actually a polynomial.  Stanley's example  is a 3-dimensional pyramid $P$ with vertices $(0,0,0)$, $(1,0,0)$, $(0,1,0)$, $(1,1,0)$, and $(1/2,0,1/2)$.  In this case, $i_P(n) = {n+3 \choose 3}$.

We say that \emph{period collapse} occurs when the minimum period is strictly less than the denominator of the polytope.  We say that $P$ has \emph{full period} if the minimum period equals the denominator of the polytope.  Stanley's example raises some natural questions.  In what dimensions can period collapse occur?  Can period collapse occur for $P$ such that $\LL(P)>2$?  What values may the minimum period be when it is not $\LL(P)$?  This note answers all of these questions.

In Section \ref{Period Collapse}, we provide (Theorem \ref{Ex}) an infinite class of 2-dimensional triangles such that, for any $\LL$, there is a triangle $P$ in this class with denominator $\LL$, but such that $i_P(n)$ is actually a polynomial.  In fact, for any $d \ge 2$ and for any $\LL$ and $s$ with $s|\LL$, there is a $d$-dimensional polytope with denominator $\LL$ but with minimum period $s$.  Such period collapse cannot occur in dimension 1, however: rational 1-dimensional polytopes always have full period (Theorem \ref{1d}).  Finally, in Section \ref{The 2-dimensional Case} (Theorem \ref{Characterize}), we give a geometric characterization of all polygons $P$ whose quasi-polynomials are actually polynomials.  We also provide several examples, one of which settles a conjecture of Zaslavsky that we detail now.

Another way to consider the period of a quasi-polynomial is to examine the periods of its coefficients.  Suppose $P$ is a $d$-dimensional polytope and, for all $j$,
\[
f_j(n)=c_{jd}n^d+c_{j,d-1}n^{d-1}+\cdots+c_{j1}n+c_{j0}.
\]
Then we say that $s_k$, the \emph{period of the $k$th coefficient}, is the minimum period of the sequence
\[
c_{1k},c_{2k},c_{3k},\ldots.
\]
The minimal period of $P$ is then the least common multiple of $s_0, s_1, \dotsc, s_d$.  T. Zaslavsky conjectured (unpublished) that the periods of the coefficients are decreasing, \textit{i.e.}, $s_{k} \leq s_{k-1}$ for $1 \leq k \leq d$.  In this paper, we provide a counterexample (Example \ref{CounterEx}) which is a 2-dimensional triangle.

\section{Period Collapse}{\label{Period Collapse}}
First, we prove that period collapse cannot happen in dimension 1.

\begin{thm}
\label{1d}
The quasi-polynomials of rational 1-dimensional polytopes always have full period.
\end{thm}

\begin{proof}In this case, $P$ is simply a segment $[\frac{p}{q}, \frac{r}{s}]$ (where the integers $p, q, r$, and $s$ are chosen so that the fractions are fully reduced).  Write $\LL = \LL(P) = \text{lcm} (s,q)$.

On the one hand, we clearly have that
\begin{equation}\label{iflooreqn}
i_P(n) = \left\lfloor n \frac{r}{s} \right\rfloor - \left\lceil n 
\frac{p}{q} \right\rceil + 1.
\end{equation}
On the other hand, there exist $\LL$ polynomials $f_1(n), \ldots, f_\LL(n)$ such that $i_P(n) = f_j(n)$, for $n \equiv j \mod{\LL}$.  The claim is that $i_P$ has period $\LL$.  To show this, it suffices to show that the constant term of $f_j(n)$ is 1 if and only if $j = \LL$.

Since $P$ is one-dimensional, we have that, for each $j\in \{1,2,\ldots,\LL\}$, the polynomial $f_j(n)$ is linear, and therefore it is determined by its values at $n=j$ and $n=j+\LL$.  Interpolating using (\ref{iflooreqn}) yields
\[
f_j(n) = \left(\frac{r}{s}-\frac{p}{q}\right)n + 1 - \left(\left\lceil j \frac{p}{q} \right\rceil - j \frac{p}{q} \right) - \left(j \frac{r}{s} - \left\lfloor j \frac{r}{s} \right\rfloor \right).
\]
The constant term is 1 if and only if $q$ and $s$ both divide $j$, which happens if and only if $j = \LL$.
\end{proof}

While in dimension 1, nothing (with respect to period collapse) is possible, in dimension 2 and higher, anything is possible, as the following theorem demonstrates.

\begin{thm}\label{Ex}
Given $d \ge 2$, and given $\LL$ and $s$ such that $s|\LL$, there exists a $d$-dimensional polytope with denominator $\LL$ whose Ehrhart quasi-polynomial has minimum period $s$.
\end{thm}

\begin{proof}
We first prove the theorem in the case where $d=2$ and $s=1$; that is, we exhibit a polygon with denominator $\LL$ for which $i_P(n)$ is actually a polynomial in $n$. Given $\LL \geq 2$, let $P$ be the triangle with vertices $(0,0), (1,\frac{\LL-1}{\LL})$, and $(\LL,0)$ (see Figure \ref{dilatedtriang}).  We will prove that
\[
i_P(n)=\frac{\LL-1}{2}n^2+\frac{\LL+1}{2}n+1.
\]

\begin{figure}
\begin{center}
\includegraphics[scale=.5]{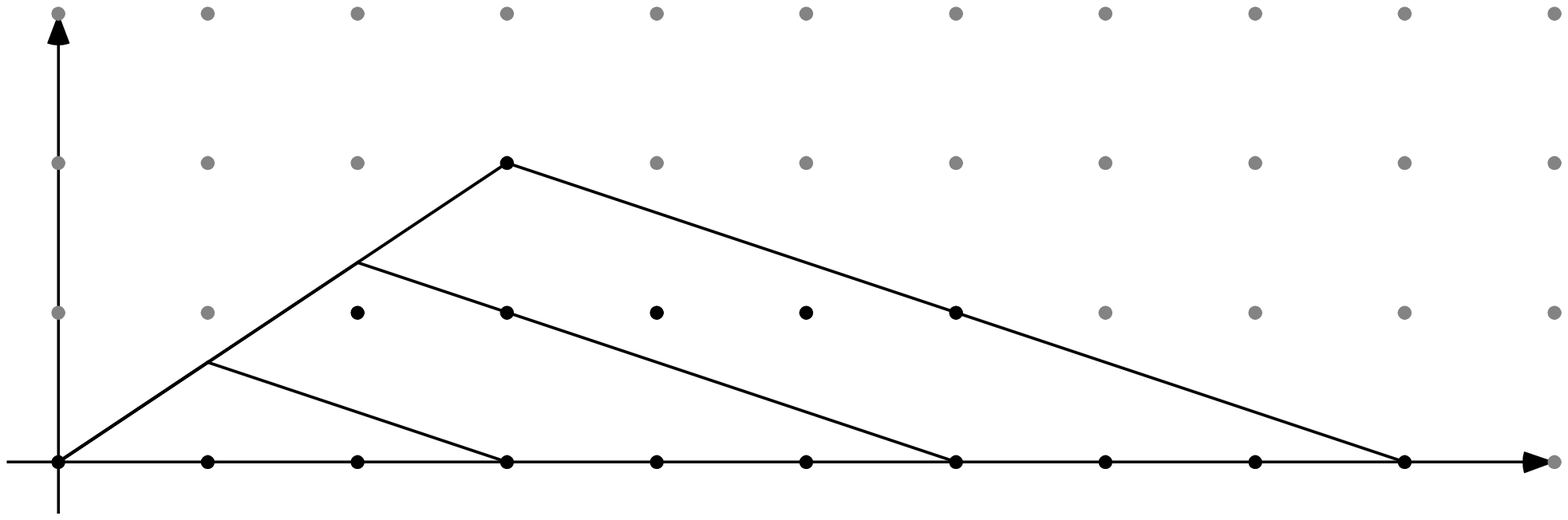}
\caption{The first three dilations of $P$ when $\LL = 3$}
\label{dilatedtriang}
\end{center}
\end{figure}

First we will calculate $i_Q(n)$, where $Q$ is the half-open parallelogram with vertices $(0,0), (1,\frac{\LL-1}{\LL}), (\LL,0),$ and $(\LL-1,-\frac{\LL-1}{\LL})$ and with top two edges open.  That is, to construct $Q$, take the closed parallelogram with these vertices and remove the line segments $\Big[(0,0), (1,\frac{\LL-1}{\LL})\Big]$ and $\Big[(1,\frac{\LL-1}{\LL}), (\LL,0)\Big]$ (see Figure \ref{dilatedparall}). $Q$ has the nice property that, for $n\in\N$, $nQ$ can be tiled by translates of $Q$ with no overlap.  It is clear that $Q$ contains exactly $\LL-1$ lattice points (the lattice points $(1,0),(2,0),\ldots,(\LL-1,0))$. To tile $nQ$, however, we must use translates of $Q$ that are not lattice translates, so it is not immediately clear how many lattice points these translates
contain. In fact, they all contain $\LL-1$ points, as we shall show.

\begin{figure}
\begin{center}
\includegraphics[scale=.5]{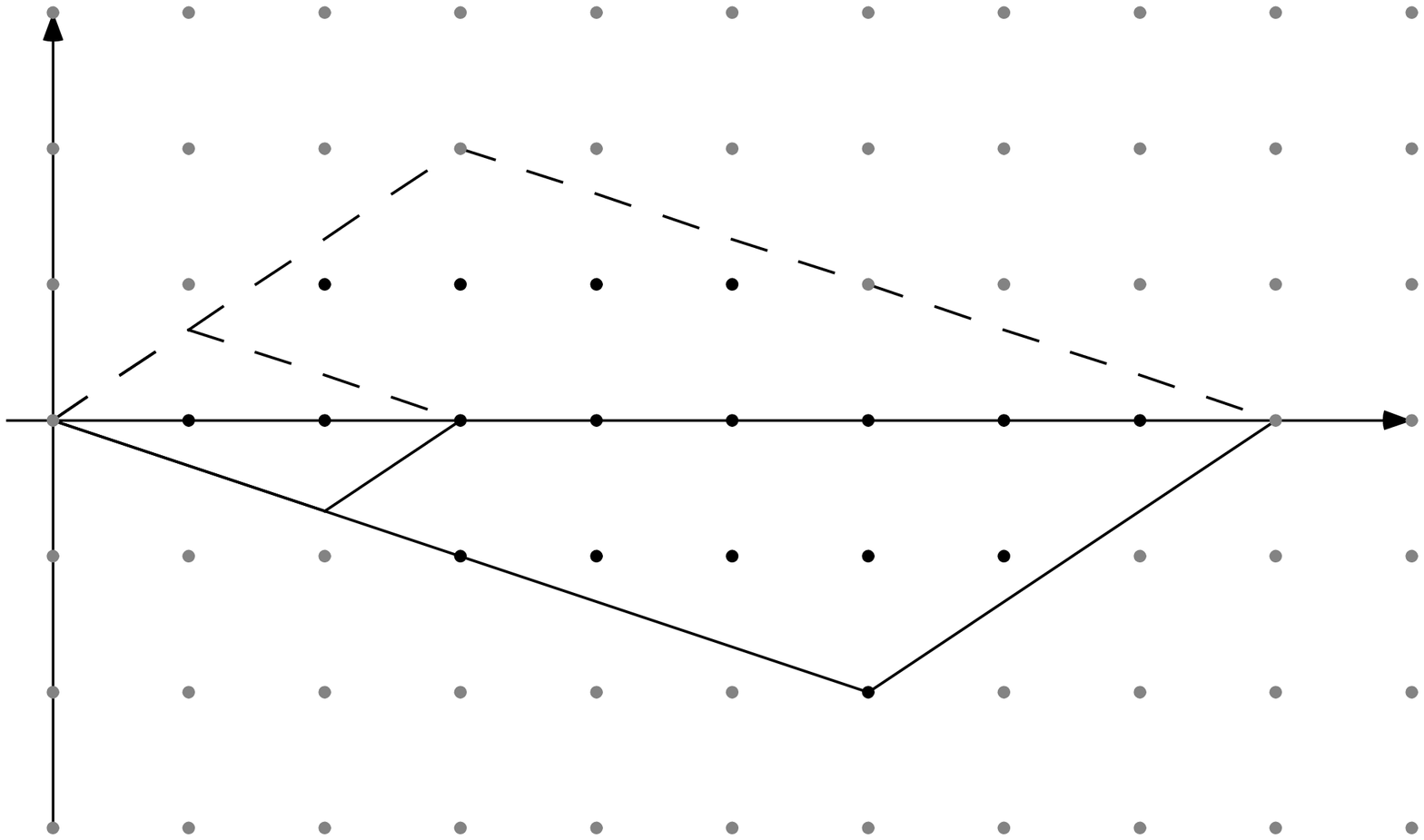}
\caption{$Q$ and $3Q$ when $\LL = 3$}
\label{dilatedparall}
\end{center}
\end{figure}

It suffices to prove this for $Q_t=Q-(0,\frac{t}{\LL})$, where $t=0,1,\ldots,\LL-1$, because all of the translates of $Q$ that we need to tile $nQ$ are lattice translates of one of these
$Q_t$. The only horizontal lines $y=a$, with $a$ integral, that possibly intersect $Q_t$ are $y=0$ and $y=-1$, and they intersect $Q_t$ with x-coordinates in the intervals $(\frac{t}{\LL-1},\LL-t)$ and $[\LL-t,\LL-1+\frac{t-1}{\LL-1}]$, respectively.  These intervals contain $\LL-t-1$ and $t$ integral points, respectively, so in all, $Q_t$ contains $\LL-1$ integer points. Therefore, we must have that \[i_Q(n)=(\LL-1)n^2.\]

Let $\overline{Q}$ be the closure of $Q$.  To calculate
$i_{\bar{Q}}(n)$, we must add to $i_Q(n)$ the number of integer
points in $n\overline{Q}\setminus nQ$, which is $n+1$ (one can
check that the number of lattice points on the interval
$\Big[(0,0),(n,n\frac{\LL-1}{\LL})\Big)$ is
$\left\lfloor{\frac{n-1}{\LL}}\right\rfloor+1$ and the number of lattice points on
the interval $\Big[(n,n\frac{\LL-1}{\LL}),(0,n\LL)\Big]$ is
$n-\left\lfloor{\frac{n-1}{\LL}}\right\rfloor$, so there are $n+1$ in all). So
\[i_{\bar{Q}}(n)=(\LL-1)n^2+n+1.\]

$n\overline{Q}$ is the union (not disjoint) of 2 copies of $nP$
(one rotated by a half-turn), each with the same number of lattice
points.  The overlap of these two copies of $nP$ is the line
segment $\Big[(0,0),(0,\LL n)\Big]$, which contains $\LL n+1$
integer points. Therefore
\[i_P(n)=\frac{1}{2}\Big(i_{\bar{Q}}(n)+(\LL n
+1)\Big)=\frac{\LL-1}{2}n^2+\frac{\LL+1}{2}n+1,\] as desired.

Now suppose $d$ is 2, but $s$ is not necessarily 1.  Let $P'$ be
the pentagon with vertices $(0,0), (1,\frac{\LL-1}{\LL}),
(\LL,0),(\LL,-\frac{1}{s}),$ and $(0,-\frac{1}{s})$.  If $P$ is
the triangle defined as before, then $nP'\setminus nP$ contains
$\left\lfloor{\frac{n}{s}}\right\rfloor\cdot(\LL n +1)$ lattice points, and so
\[i_{P'}(n)=i_P(n)+\left\lfloor{\frac{n}{s}}\right\rfloor\cdot(\LL n
+1),\] which has minimum period $s$.

Now suppose $d$ is greater than 2.  Let $P'$ be the pentagon
defined as before, and let $P''=P'\times[0,1]^{d-2}$, a polytope
of dimension $d$.  Then
\[i_{P''}(n)=(n+1)^{d-2}i_{P'}(n),\]
which also has minimum period $s$.
\end{proof}

\section{The 2-dimensional Case}{\label{The 2-dimensional Case}}
We have seen (in Theorem \ref{Ex}) an infinite class of rational polygons $P$ in dimension 2 such that $i_P(n)$ is a polynomial.  Can we characterize such polygons?  We know that, for all \emph{integer} polygons $P$,  $i_P(n)$ is a polynomial.  One
property that an integer polygon $P$ has is that it and its dilates
satisfy Pick's theorem, i.e., if we let $\partial_P(n)=\#\Big(boundary(nP)\cap \Z^d\Big)$, then
\begin{align*}
i_P(n) &= \text{Area}(nP)+\frac{1}{2}\partial_P(n)+1\\
&= n^2\text{Area}(P)+\frac{1}{2}\partial_P(n)+1.
\end{align*}
Another
property that an integer polygon, $P$, and its dilates satisfy is that the number of
points on their boundary is linear, i.e.,
\[\partial_P(n)=n\partial_P(1).\]
In fact, these two properties are exactly what we need to
guarantee that a rational polygon's Ehrhart quasi-polynomial is
actually a polynomial.

\begin{thm}
\label{Characterize}
Let $P\subset \Z^2$ be a rational polygon, let $A$ be the area of
$P$, and let $\LL$ be the denominator of $P$.  Then the following are equivalent:
\begin{enumerate}
    \item $i_P(n)$ is a polynomial in n;
    \item $i_P(n)=An^2+\frac{1}{2}\partial_P(1)n+1$;
    \item For all $n\in\N$,
    \begin{enumerate}
        \item nP obeys Pick's theorem, i.e.,
        $i_P(n)=An^2+\frac{1}{2}\partial_P(n)+1$, and
        \item $\partial_P(n)=n\partial_P(1)$; and
    \end{enumerate}
    \item For $n=1,2,\ldots,\LL$, 3a and 3b hold.
\end{enumerate}
\end{thm}
\begin{proof}
We will prove that $1\Rightarrow 2 \Rightarrow 3 \Rightarrow 4
\Rightarrow 2 \Rightarrow 1$.  Two of these steps, $3 \Rightarrow
4$ and $2 \Rightarrow 1$, are trivial.  To prove the remaining implications, we will repeatedly use the law of reciprocity for Ehrhart quasi-polynomials, which was stated in the introduction.

$1\Rightarrow 2$.  If 1 holds, then $i_P(n)=An^2+bn+c$ for some
$b$ and $c$.  Since $i_P(0)=1$, we know that $c=1$.
By the reciprocity law, we know that
\[i^\circ_P(n)=A(-n)^2+b(-n)+c,\]
and so
\[\partial_P(1)=i_P(1)-i^\circ_P(1)=2b.\]
Therefore $i_P(n)=An^2+\frac{1}{2}\partial_P(1)n+1$, as desired.

$2\Rightarrow 3$.  If 2 holds, then, again using reciprocity, for
all $n\in \N$,
\[i^\circ_P(n)=An^2-\frac{1}{2}\partial_P(1)n+1,\]
and so
\[\partial_P(n)=i_P(n)-i^\circ_P(n)=\partial_P(1)n,\]
and so 3b holds.  Then
\begin{align*}
i_P(n)&=An^2+\frac{1}{2}\partial_P(1)n+1\\
&=An^2+\frac{1}{2}\partial_P(n)+1,
\end{align*}
and so 3a holds.

$4\Rightarrow 2$.  If 4 holds, then let
\[f_j(n)=An^2+b_jn+c_j,\]
for $j=1,2,\ldots,\LL$, be the polynomials such that
$i_P(n)=f_j(n)$ for $n \equiv j \mod{\LL}$.  Given $j$ with $1\le
j\le \LL$, we again use reciprocity, and we have
\begin{equation}\label{a}
\begin{split}
j\partial_P(1)&=\partial_P(j)\\
&=f_j(j)-f_{\LL-j}(-j)\\
&=(b_j+b_{\LL-j})\cdot j+(c_j-c_{\LL-j})
\end{split}
\end{equation}
and
\begin{equation}\label{b}
\begin{split}
(\LL-j)\partial_P(1)&=\partial_P(\LL-j)\\
&=f_{\LL-j}(\LL-j)-f_{j}(j-\LL)\\
&=(b_j+b_{\LL-j})\cdot (\LL-j)+(c_{\LL-j}-c_j)
\end{split}
\end{equation}
Multiplying Equation \eqref{a} by $\LL-j$ and Equation \eqref{b}
by $j$ and subtracting,
\[0=\LL\cdot (c_j-c_{\LL-j}),\]
and so
\begin{equation}\label{e}
c_j=c_{\LL-j}.
\end{equation}
Adding Equations \eqref{a} and \eqref{b},
\[\LL\cdot \partial_P(1)=\LL\cdot(b_j+b_{\LL-j}),\]
and so
\begin{equation}\label{f}
b_j+b_{\LL-j}=\partial_P(1).
\end{equation} Using the facts that
Pick's theorem holds and that $j\partial_P(1)=\partial_P(j)$, we
have
\begin{align*}
Aj^2+\frac{1}{2}\partial_P(1)\cdot j+1&=Aj^2+\frac{1}{2}\partial_P(j)+1\\
&=f_j(j)\\
&=Aj^2+b_j\cdot j+c_j,
\end{align*}
and so
\begin{equation}\label{c}
\frac{1}{2}\partial_P(1)\cdot j+1=b_j\cdot j+c_j.
\end{equation}
Similarly,
\begin{equation}\label{d}
\frac{1}{2}\partial_P(1)\cdot(\LL-j)+1=b_{\LL -
j}\cdot(\LL-j)+c_{\LL-j}.
\end{equation}
Multiplying Equation \eqref{c} by $\LL-j$ and Equation \eqref{d}
by $j$ and adding together (and then using Equations \eqref{e} and
\eqref{f}),
\begin{align*}
\partial_P(1)\cdot j \cdot (\LL-j) + \LL &=
(b_j+b_{\LL-j})\cdot j\cdot (\LL-j) +(\LL-j)\cdot c_j+j\cdot c_{\LL-j}\\
&=\partial_P(1)\cdot j \cdot (\LL-j) + \LL\cdot c_j,
\end{align*}
and so $c_j=1$.  Substituting $c_j=1$ into Equation \eqref{c}, we see
that $b_j=\frac{1}{2}\partial_P(1)$.  Therefore, for all $n\in
\N$,
\[i_P(n)=An^2+\frac{1}{2}\partial_P(1)n+1,\]
as desired.
\end{proof}

\begin{exm}
$P$ is the triangle with vertices $(0,0), (\LL,0),$ and
$(1,\frac{\LL-1}{\LL})$, for some $\LL\in \N$.
\end{exm}
This is the example from Theorem \ref{Ex} with denominator $\LL$ for which the Ehrhart
quasi-polynomial is a polynomial. One can check that conditions 3a and 3b are met.

\begin{exm}
\label{CounterEx}
$P$ is the triangle with vertices $(-\frac{1}{2},-\frac{1}{2}),
(\frac{1}{2},-\frac{1}{2}),$ and $(0,\frac{3}{2})$.
\end{exm} One
can check that $nP$, for $n\in \N$, satisfies 3a (Pick's theorem),
but not 3b. Indeed, we have
\[i_P(n)=\left\{%
\begin{array}{ll}
    n^2+1, & \hbox{if $n$ is odd} \\
    n^2+n+1, & \hbox{if $n$ is even,} \\
\end{array}%
\right.\] which is not a polynomial. This example disproves a
conjecture of T. Zaslavsky that the period of the coefficient of
$n^k$ in the quasi-polynomial increases as $k$ decreases (in the
example, the coefficients of $n^2$ and $n^0$ have period 1, but
the coefficient of $n^1$ has period 2).  A similar counterexample
has been found independently by D. Einstein.

\begin{exm}
$P$ is the triangle with vertices $(0,0),(1,0),$ and
$(0,\frac{1}{2})$.
\end{exm}
In this example, $nP$, for $n\in \N$ satisfies 3b, but not 3a.  We
have
\[i_P(n)=\left\{%
\begin{array}{ll}
    \frac{1}{4}n^2+n+\frac{3}{4}, & \hbox{if $n$ is odd,} \\
    \frac{1}{4}n^2+n+1, & \hbox{if $n$ is even.} \\
\end{array}%
\right.\]

\section*{Acknowledgements}
We would like to thank Matthias Beck and Jesus De Loera for helpful conversations.  Special thanks to David Einstein for a simplification of the example in Theorem \ref{Ex}.


\begin{thebibliography}{99}

\bibitem{ehrhart} \textsc{E. Ehrhart}, \emph{Polyn™mes arithmŽtiques et mŽthode des polydres en combinatoire}. International Series of Numerical Mathematics, Vol. 35. BirkhŠuser Verlag, Basel-Stuttgart, 1977.

\bibitem{macdonald} \textsc{I.\,G. MacDonald}, Polynomials associated with finite cell complexes,  \emph{J. London Math. Soc.}  \textbf{4} (1971), 181-192.

\bibitem{stanley} \textsc{R.\,P. Stanley}, \emph{Enumerative Combinatorics, Volume I}, Cambridge University Press, 1997.

\end{thebibliography}
\end{document}